\documentclass[a4paper,12pt]{article}
\usepackage{amsmath}
\usepackage{amsthm}
\usepackage{amssymb}
\usepackage{amscd}
\usepackage{verbatim}
\usepackage{latexsym}
\theoremstyle{definition}
\numberwithin{equation}{section}

\newtheorem{lmm}{Lemma}[section]

\newtheorem{prp}{Proposition}[section]

\newtheorem{thm}{Theorem}[section]

\theoremstyle{definition}
\newtheorem{rem}{Remark}
\newtheorem{dfn}{Definition}[section]
\newtheorem{exa}{Example}[section]
\theoremstyle{remark}
\newtheorem{note}{Note}

\def\dim{\mathop{\mathrm{dim}}\nolimits}

\def\Ind{\mathop{\mathrm{Ind}}\nolimits}
\def\Res{\mathop{\mathrm{Res}}\nolimits}

\newcommand{\mf}[1]{{\mathfrak{#1}}}

\newcommand{\bb}[1]{{\mathbb{#1 \ }}}
\title{\textbf{\sf{Classification of the Irreducible Representations of the Affine Hecke algebra of Type $B_2$ with Unequal Parameters}}}
\author{Naoya Enomoto\\
Research Institute for Mathematical Science\\
Kyoto University \\
henon@kurims.kyoto-u.ac.jp}
\date{}
\begin{document}
\maketitle 
\quad {} \\
\section{Introduction}
 \indent  The representation theory of the affine Hecke algebras has two different approaches. One is a geometric approach and the other is a combinatorial one. \\
\indent In the equal parameter case, affine Hecke algebras are constructed using equivariant K-groups, and their irreducible representations are constructed on Borel-Moore homologies. By this method, their irreducible representations are parameterized by the index triples (\cite{CG},\cite{KL}). On the other hand, G. Lusztig classified the irreducible representations in the unequal parameter case. His ideas are to use equivariant cohomologies and graded Hecke algebras (\cite{L1989},\cite{L1},\cite{L2},\cite{L3}). \\
\indent Although the geometric approach will give us a powerful method for the classification, but it does not tell us the detailed structure of irreducible representations. Thus it is important to construct them explicitly in combinatorial approach. \\
\indent  Using semi-normal representations and the generalized Young tableaux, A. Ram constructed calibrated irreducible representations with equal parameters (\cite{R1}). Furthermore C. Kriloff and A. Ram constructed irreducible calibrated representations of graded Hecke algebras (\cite{KR}). However, we cannot always construct irreducible representations by combinatorial manner.\\
\indent A. Ram classified irreducible representations of affine Hecke algebras of type $A_1$, $A_2$, $B_2$, $G_2$ in equal parameter case (\cite{R2}). But there are some mistakes in his list of irreducible representations and his construction of induced representation of type $B_2$.\\
\indent In this paper, we will correct his list about type $B_2$ and also classify the irreducible representations in the unequal parameter case. There are three one-parameter families of calibrated irreducible representations  and some other irreducible representations. 
\\
\textbf{Acknowledgement}. I would like to thank Professor M. Kashiwara and Professor S. Ariki for their advices and suggestions, and Mathematica for its power of calculation.\\
\section{Preliminaries}
\subsection{Affine Hecke algebra}
 We will use following notations.
\[
\begin{array}{ll}
(R,R^{+},\Pi) & \text{a root system of finite type, its positive roots and simple roots},\\
Q,P & \text{the root lattice and the weight lattice of R},\\
Q^{\vee},P^{\vee} & \text{the coroot lattice and the coweight lattice of R} \\ 
W & \text{the Weyl group of R}, \\
\ell(w) & \text{the length of $w \in W$}
\end{array}
\]
We put $\Pi=\{\alpha_i\}_{i \in I}$, and denote by $s_i$ the simple reflection associated with $\alpha_i$.
\quad {} \\
First we define the Iwahori-Hecke algebra of W.
\begin{dfn}
Let $\{q_i\}_{i \in I}$ be indeterminates. Then the \textit{Iwahori-Hecke algebra} $\mathcal{H}$ of W is the associative algebra over $\bb{C}(q_i)$ defined by following generators and relations;
\begin{eqnarray*}
\text{generators} \ &{}& T_i \quad (i \in I) \\
\text{relations} \ &{}& (T_i-q_i)(T_i+q_i^{-1})=0 \quad (i \in I), \\
\quad &{}& \overbrace{T_iT_jT_i \cdots}^{m_{ij}}=\overbrace{T_jT_iT_j \cdots}^{m_{ij}},
\end{eqnarray*}
where $m_{ij}=2,3,4,6$ according to $\langle \alpha_i,\alpha_j^{\vee}\rangle \langle \alpha_j \alpha_i^{\vee}\rangle=0,1,2,3$.
\begin{rem}
Indeterminates $q_i,q_j$ must be equal if and only if $\alpha_i,\alpha_j$ are in the same $W$-orbit in $R$. If all $q_i$ are equal, we call the \textit{equal parameter case}, and otherwise, \textit{ the unequal parameter case}.
\end{rem}
\end{dfn}
\indent  For a reduced expression $s_{i_1}s_{i_2} \cdots s_{i_r}$ of $w \in W$, we define $T_w=T_{i_1}T_{i_2} \cdots T_{i_r}$.  This does not depend on the choice of reduced expressions.\\
\indent Let us define the affine Hecke algebras.
\begin{dfn}
The \textit{affine Hecke algebra} $\widehat{\mathcal{H}}$ is the associative algebra over $\bb{C}(q_i;i \in I)$ defined by following generators and relations;
\begin{eqnarray*}
\text{generators} \  &{}& T_wX^\lambda \quad (w \in W,\lambda \in P^\vee), \\
\text{relations} \ &{}& (T_i-q_i)(T_i+q_i^{-1})=0 \quad (i \in I), \\
\quad &{}& T_wT_{w'}=T_{ww'} \quad \text{if} \quad \ell(w)+\ell(w')=\ell(ww') \ (w,w' \in W), \\
\quad &{}& X^\lambda{X^\mu}=X^{\lambda+\mu} \quad (\lambda,\mu \in P^{\vee}), \\
\quad &{}& X^\lambda{T_i}=T_iX^{s_i\lambda}+(q_i-q_i^{-1})\frac{X^\lambda-X^{s_i\lambda}}{1-X^{-\alpha_i^\vee}} \quad (i \in I).
\end{eqnarray*}
\end{dfn}

\subsection{Principal series representations and their irreducibility}
Let us put $X^{P^\vee}=\{X^\lambda|\lambda \in P^{\vee}\}$ and let $\chi:X^{P^\vee} \to \bb{C}^*$ be a character of $X^{P^{\vee}}$.
\begin{dfn}
Let $\bb{C}v_\chi$ be the one-dimensional representation of $\bb{C}[X]$ defined by
\[
X^\lambda \cdot v_\chi=\chi(X^\lambda)v_\chi.
\]
We call $M(\chi)=\Ind_{\bb{C}[X]}^{\widehat{\mathcal{H}}}\bb{C} v_\chi=\widehat{\mathcal{H}} \otimes_{\bb{C}[X]} \bb{C}v_\chi$ the \textit{principal representation of $\mathcal{\widehat{H}}$} associated with $\chi$. 
\end{dfn}
Note that $\Res_{\mathcal{H}}^{\widehat{\mathcal{H}}}M(\chi)$ is isomorphic to the regular representation of $\mathcal{H}$, so that $\dim{M(\chi)}=|W|$. \\
\indent We put
\[
q_\alpha=q_i\ \text{for} \ \alpha^{\vee} \in W\alpha_i^{\vee} \ (i \in I).
\]
\begin{thm}[Kato's Criterion of Irreducibility]\label{kato}
Let us put
\[
P(\chi)=\{\alpha^{\vee}>0|\chi(X^{\alpha^\vee})=q_\alpha^{\pm{2}}\}.
\]
Then $M(\chi)$ is irreducible if and only if $P(\chi)=\phi$.
\end{thm}
For any finite-dimensional representation of $\widehat{\mathcal{H}}$ we put
\begin{eqnarray*}
M_\chi&=&\{v \in M|X^\lambda{v}=\chi(X^\lambda)v \ \text{for any}X^\lambda \in X\}, \\
M_\chi^{\text{gen}}&=&\left\{v \in M\left|
\begin{array}{l}
\text{there exists} \ k>0 \ \text{such that} \\
(X^\lambda-\chi(X^\lambda))^kv=0 \ \text{for any} \  X^\lambda \in X
\end{array}
\right.\right\}.
\end{eqnarray*}
Then $\displaystyle M=\bigoplus_{\chi \in T}M_t^{\text{gen}}$ is the generalized weight decomposition of M. \\
\begin{prp}
If $M$ is a simple $\widehat{\mathcal{H}}$-module with $M_\chi \neq 0$, then $M$ is a quotient of $M(\chi)$.
\end{prp}
\begin{dfn}
A finite-dimensional representation $M$ of $\widehat{\mathcal{H}}$ is \textit{calibrated}(or \textit{X-semisimple}) if $M_\chi^{\text{gen}}=M_\chi \ (\text{for all} \ \chi)$.
\end{dfn}
\subsection{W-action Lemma}
Let us define the action of Weyl group $W$ as the following;
\[
(w \cdot \chi)(X^\lambda)=\chi(X^{w^{-1}\lambda}) \ (w \in W,\lambda \in P^{\vee}).
\]
The following proposition is well known.
\begin{prp}[W-action Lemma]\label{WW} \quad {} \\
(1) \ If $M(\chi) \cong M(\chi')$, then there exists $w \in W$ such that $\chi=w\chi$. \\
(2) \ The representations $M(\chi)$ and $M(w\chi)$ have the same composition factors. 
\end{prp}
\subsection{Specialization lemma}
\indent Let $\bb{K}$ be a field and $\bb{S}$ a discrete valuation ring such that $\bb{K}$ is the fraction field of $\bb{S}$. Let us denote the $\mf{m}=(\pi)$ the maximal ideal of $\bb{S}$ and let $\bb{F}=\bb{S}/\mf{m}$ be the residue field of $\bb{S}$. Let $K(\widehat{\mathcal{H}}_{\bb{F}}\text{-mod})$ be the Grothendieck group of the category of finite-dimensional representations of $\widehat{\mathcal{H}}_{\bb{F}}$. \\
\indent the following lemma is well-known (e.g. see \cite[Lemma 13.16]{A}.)
\begin{lmm}[Specialization Lemma]\label{spl}
\ Let $V$ be an $\widehat{\mathcal{H}}_{\bb{K}}$-module and $L$ an $\mathcal{\widehat{H}}_{\bb{S}}$-submodule of $V$ which is an $\bb{S}$-lattice of full rank. Then $[L \otimes \bb{F}] \in K(\widehat{\mathcal{H}}_{\bb{F}}\text{-mod})$ is determined by $V$ and does not depend on the choice of $L$.
\end{lmm}
\subsection{Key results for type $B_2$}
 Let us consider the type $B_2$;
\begin{eqnarray*}
&{}&P^{\vee}=\bb{Z}\varepsilon_1 \oplus \bb{Z}\varepsilon_2, \ R^{\vee}=\{\alpha_1^{\vee}=\varepsilon_1-\varepsilon_2,\alpha_2^{\vee}=2\varepsilon_2\}, \ X_i=X^{\varepsilon_i}.\\
&{}&s_1\varepsilon_1=\varepsilon_2, \ s_1\varepsilon_2=\varepsilon_1, \quad s_2\varepsilon_1=\varepsilon_1, \ s_2\varepsilon_2=-\varepsilon_2
\end{eqnarray*}
Let us recall the definition of affine Hecke algebra of type $B_2$ with unequal parameters.
\begin{dfn}
The \textit{affine Hecke algebra$\widehat{\mathcal{H}}$ of type $B_2$} is the associative algebra over $\bb{C}(p,q)$ defined by the following generators and relations;
\begin{eqnarray*}
\begin{array}{lll}
\text{generators} \  & T_1,T_2,X_1,X_2 & \\
\text{relations} \ & (T_1-q)(T_1+q^{-1})=0, & (T_2-p)(T_2+p^{-1})=0, \\
 & T_1T_2T_1T_2=T_2T_1T_2T_1, & \\
& T_1X_2T_1=X_1, & T_2X_2^{-1}T_2=X_2, \\
& T_2X_1=X_1T_2, & X_1X_2=X_2X_1. 
\end{array}
\end{eqnarray*} 
\end{dfn}
\indent We will use the following four subalgebras of $\widehat{\mathcal{H}}(B_2)$;
\begin{eqnarray*}
\widehat{\mathcal{H}}_1=\langle{T_1,X_1,X_2}\rangle, \quad 
\widehat{\mathcal{H}}_2=\langle{T_2,X_1,X_2}\rangle, \quad 
\mathcal{H}=\langle{T_1,T_2}\rangle, \quad 
\bb{C}[X_1,X_2]\subset \widehat{\mathcal{H}}.
\end{eqnarray*}
\begin{lmm}[Decomposition Lemma]\label{dec}
Suppose $\chi(X^{\alpha_i})=q_i^2$, and let $\rho_1,\rho_2$ be the following 1-dimensional representations of $\widehat{\mathcal{H}}_i=\langle T_i,X_j(1 \le j \le 2) \rangle \subset \widehat{\mathcal{H}}$;
\[
\rho_1(X_j)=\chi(X_j), \ \rho_1(T_i)=q_i, \ \rho_2(X_j)=(s_i\chi)(X_j), \ \rho_2(T_i)=-q_i^{-1}.
\]
Then there exists the following short exact sequence;
\[
0 \to \Ind_{\widehat{\mathcal{H}}_i}^{\widehat{\mathcal{H}}}\rho_2 \to M(\chi) \to \Ind_{\widehat{\mathcal{H}}_i}^{\widehat{\mathcal{H}}}\rho_1 \to 0
\] 
\end{lmm}
\section{Classification}
\subsection{Method}
\indent Let $M$ be an irreducible representation which is not principal. Then $M$ appears in some $M(\chi)$. By Kato's criterion (Theorem \ref{kato}), $P(\chi) \neq \phi$. Using W-action Lemma (Lemma \ref{WW}), we may assume $P(\chi) \ni \alpha_1$ or $\alpha_2$. thus we obtain the following Lemma. We will use the notation $-\chi$ defined by $(-\chi)(X_i)=-\chi(X_i) \ (i=1,2)$.
\begin{lmm}
Except irreducible principal series representations, any finite-dimensional irreducible representation appears in the principal representations associated with the following characters as their composition factors;
\[
\begin{array}{|c||c|c|c|c|c|c|c|c|c|c|}
\hline
\chi & \chi_a & \chi_b & \chi_c & \chi_d^{(1)} & \chi_d^{(2)} & \chi_d^{(3)} & \chi_d^{(4)} & \chi_d^{(5)} & \chi_f(v) & \chi_g(u) \\
\hline
\hline
\chi(X_1) & q^2p & q^2p^{-1} & -p^{-1} & q^2 & q & p & 1 & 1 & pv & q^2u \\
\hline
\chi(X_2) & p & p^{-1} & p & 1 & q^{-1} & p & p & p & p & u \\
\hline
\end{array}
\] 
and $-\chi_a,-\chi_b,-\chi_d^{(1)},-\chi_d^{(2)},-\chi_d^{(3)},-\chi_d^{(4)},-\chi_d^{(5)},-\chi_f(v)$, where 
\begin{eqnarray*}
v &\neq& \pm{p^{-2}},\pm{p^{-1}},\pm{1},q^{\pm{2}},q^{\pm{2}}p^{-2}, \\
u &\neq& \pm{p^{\pm{1}}},\pm{1},\pm{q^{-2}},\pm{q^{-1}},\pm{q^{-2}}p^{\pm{1}}.
\end{eqnarray*}
\end{lmm}
\begin{note}
Two principal series representations $M(-\chi_c)$ and $M(\chi_c)$ have same composition factors, because of $W$-action lemma (Lemma \ref{WW}). By replacing $u$ with $-u$, we don't need to consider  $-\chi_g(u)$.
\end{note}
\indent Finaly,  we must determine the composition factors of $M(\chi)$ for above characters, and we must prove their irreducibility. But using the decomposition lemma, we consider the representations induced from $\mathcal{\widehat{H}}_i$. We will show the examples and some proofs in the following section. 
\subsection{Some examples and proofs}
\begin{exa}\label{exd5}
We consider the principal series representation $M(\chi_d^{(5)})$. Let $\rho_1^{d^{(5)}}$ and $\rho_2^{d^{(5)}}$ be the following 1-dimensional representations of $\widehat{\mathcal{H}}_2$;
\[
\begin{array}{|c||c|c|c|}
\hline
 & X_1 & X_2 & T_2 \\
\hline
\hline
\rho_1^{d^{(5)}} & -1 & p & p \\
\hline
\rho_2^{d^{(5)}} & -1 & -p^{-1} & -p^{-1}\\
\hline
\end{array}
\]
Since $\chi_d^{(5)}(\alpha_2^{\vee})=p^2$, we can apply the decompose lemma (Lemma \ref{dec}) to $M(\chi_d^{(5)})$.
\begin{lmm}\label{lmd5}
Suppose $p \neq -q^{\pm{2}}$. Then $\Ind_{\widehat{\mathcal{H}}_2}^{\widehat{\mathcal{H}}}\rho_1^{d^{(5)}}$ and $\Ind_{\widehat{\mathcal{H}}_2}^{\widehat{\mathcal{H}}}\rho_2^{d^{(5)}}$ are 4-dimensional non-calibrated irreducible representations.  
\end{lmm}
\begin{proof}
We consider the case of $\Ind_{\widehat{\mathcal{H}}_2}^{\widehat{\mathcal{H}}}\rho_1^{d^{(5)}}$. These simultaneous eigenvalues of $X_1$ and $X_2$ are $(p,-1),(-1,p)$, and the multiplicity of each eigenvalues is two. We can find the following representation matrices;
\begin{eqnarray*}
T_1&=&\left(
\begin{smallmatrix}
\frac{p(q^2-1)}{(1+p)q} & -\frac{(p-1)(q^2-1)}{(1+p)q}& 1 & -\frac{p(q^2-1)^2}{(1+p)^2q^2}\\
0 & \frac{p(q^2-1)}{(1+p)q}& 0 & \frac{(p+q^2)(1+pq^2)}{(1+p)^2q^2}\\
\frac{(p+q^2)(1+pq^2)}{(1+p)^2q^2} & \frac{(1-p+p^2)(q^2-1)^2}{(1+p)^2q^2}& \frac{(q^2-1)}{(1+p)q}& \frac{(p-1)(q^2-1)(p+q^2)(1+pq^2)}{(1+p)^3q^3}\\
0 & 1 & 0 & \frac{(q^2-1)}{(1+p)q} 
\end{smallmatrix}
\right), \\
T_2&=&\left(
\begin{smallmatrix}
 & 1 & & \\
1 & \frac{(p^2-1)}{p} & & \\
 & & p & \\
 & & & p
\end{smallmatrix}
\right),\\
X_1&=&\left(
\begin{smallmatrix}
p & & & \\
 & p & & \\
& & -1 & -\frac{(-1+p)(p+q^2)(1+pq^2)}{p(1+p)q^2}\\
& & & -1
\end{smallmatrix}
\right),\ 
X_2=\left(
\begin{smallmatrix}
-1 & -\frac{p^2-1}{p}& & \\
 & -1 & & \\
& & p & \\
& & & p
\end{smallmatrix}
\right).
\end{eqnarray*}
Since $p \neq -q^{\pm{2}}$ and $p,q$ are not a root of unity, the non-diagonal component with respect to $(p,-1),(-1,p)$ in $X_1$ and $X_2$ don't vanish. Thus the dimension of each simultaneous eigenspaces is just one. Let $v_1,v_2$ be the simultaneous eigenvectors with respect to $(p,-1),(-1,p)$. We have
\begin{eqnarray*}
T_1v_1=\frac{p(q^2-1)}{(1+p)q}v_1+\frac{(p+q^2)(1+pq^2)}{(1+p)^2q^2}v_2, \quad
T_1v_2=\frac{q^2-1}{(1+p)q}v_2+v_1,
\end{eqnarray*} 
and $p \neq -q^{\pm{2}}$.  If there exists a submodule $0 \neq U$ of $\Ind_{\widehat{\mathcal{H}}_2}^{\widehat{\mathcal{H}}}\rho_1^{d^{(5)}}$, then $U$ contains $v_1$ or $v_2$. If $v_2$ is conteined in $U$, then $v_1$ is contained in $U$, and vice versa. Therefore $\langle v_1,v_2,T_2v_1,T_1T_2v_1 \rangle \subset U$. This implies that $U=\Ind_{\widehat{\mathcal{H}}_2}^{\widehat{\mathcal{H}}}\rho_1^{d^{(5)}}$, and $\Ind_{\widehat{\mathcal{H}}_2}^{\widehat{\mathcal{H}}}\rho_1^{d^{(5)}}$ is irreducible. Similarly, we can show that $\Ind_{\widehat{\mathcal{H}}_2}^{\widehat{\mathcal{H}}}\rho_2^{d^{(5)}}$ is irreducible.
\end{proof}
\end{exa}

\begin{exa}\label{exa}
We consider $M(\chi_a)$. Let $\rho_1^{a}$ and $\rho_2^{a}$ be the following 1-dimensional representations of $\widehat{\mathcal{H}}_2$;
\[
\begin{array}{|c||c|c|c|}
\hline
 & X_1 & X_2 & T_2 \\
\hline
\hline
\rho_1^{a} & q^2p & p & p \\
\hline
\rho_2^{a} & q^2p & -p^{-1} & -p^{-1}\\
\hline
\end{array}
\]
Since $\chi_a(\alpha_2^{\vee})=p^2$, we can apply the decompose lemma (Lemma \ref{dec}) to $M(\chi_a)$. 
\begin{lmm}\label{lma}
Suppose $p \neq \pm{q^{-1}},\pm{q^{-2}},p^2 \neq -q^{-2}$. Then $\Ind_{\widehat{\mathcal{H}}_2}^{\widehat{\mathcal{H}}}\rho_1^{a}$ and $\Ind_{\widehat{\mathcal{H}}_2}^{\widehat{\mathcal{H}}}\rho_2^{a}$ have 1- and 3-dimensional calibrated irreducible composition factors. More precisely,\\
(1) \ $\Ind_{\widehat{\mathcal{H}}_2}^{\widehat{\mathcal{H}}}\rho_1^{a}$ have two composition factors which are presented by the following representation matricies;
\begin{itemize}
\item{$X_1=pq^2,\ X_2=p, \ T_1=q, \ T_2= p$.}
\item{$U_a^1$:\begin{eqnarray*}
&{}&X_1=\left(\begin{smallmatrix}
p & & \\
& p & \\
 & & p^{-1}q^{-2}
\end{smallmatrix}\right), \ 
X_2=\left(\begin{smallmatrix}
pq^2 & & \\
& p^{-1}q^{-2} & \\
 & & p
\end{smallmatrix}\right), \\
&{}&T_1=\left(\begin{smallmatrix}
-q^{-1} & & \\
 & \frac{p^2q(q^2-1)}{(p^2q^2-1)}& \frac{(p^2-1)(p^2q^4-1)}{(p^2q^2-1)^2}\\
 & 1 & -\frac{(q^2-1)}{q(p^2q^2-1)} 
\end{smallmatrix}\right),\ T_2=\left(\begin{smallmatrix}
\frac{p(p^2-1)q^4}{(p^2q^4-1)} & \frac{(q^4-1)(p^4q^4-1)}{(p^2q^4-1)^2}& \\
1 & -\frac{p^2-1}{p(p^2q^4-1)}& \\
 & & p
\end{smallmatrix}\right).
\end{eqnarray*}}
\end{itemize}
(2) \ $\Ind_{\widehat{\mathcal{H}}_2}^{\widehat{\mathcal{H}}}\rho_2^{a}$ have two composition factors which are presented by the following representation matricies;
\begin{itemize}
\item{$X_1=p^{-1}q^{-2}, \ X_2=p^{-1}, \ T_1=-q^{-1}, \ T_2=-p^{-1}$.}
\item{$U_a^2$:
\begin{eqnarray*}
&{}&X_1=\left(\begin{smallmatrix}
p^{-1} & & \\
& pq^2 & \\
 & & p^{-1}
\end{smallmatrix}\right), \ X_2=\left(\begin{smallmatrix}
pq^2 & & \\
& p^{-1} & \\
 & & p^{-1}q^{-2}
\end{smallmatrix}\right), \\
&{}&T_1=\left(\begin{smallmatrix}
-\frac{q^2-1}{q(p^2q^2-1)} & 1& \\
\frac{(p^2-1)(p^2q^4-1)}{(p^2q^2-1)^2} & \frac{p^2q(q^2-1)}{(p^2q^2-1)}& \\
 & & q
\end{smallmatrix}\right), \ T_2=\left(\begin{smallmatrix}
\frac{p(p^2-1)}{(p^2q^4-1)} & & \frac{(q^4-1)(p^4q^4-1)}{(p^2q^4-1)^2}\\
 & -p^{-1} & \\
1 & & -\frac{p^2-1}{p(p^2q^4-1)}
\end{smallmatrix}\right).
\end{eqnarray*}
}
\end{itemize}
\end{lmm}
\end{exa}
\begin{exa}\label{exb}
We consider $M(\chi_b)$. Let $\rho_1^{b}$ and $\rho_2^{b}$ be the following 1-dimensional representations of $\widehat{\mathcal{H}}_1$;
\[
\begin{array}{|c||c|c|c|}
\hline
 & X_1 & X_2 & T_1 \\
\hline
\hline
\rho_1^{b} & q^2p^{-1} & p^{-1} & q \\
\hline
\rho_2^{b} & p^{-1} & q^2p^{-1} & -q^{-1}\\
\hline
\end{array}
\]
Since $\chi_a(\alpha_1^{\vee})=q^2$, we can apply the decompose lemma (Lemma \ref{dec}) to $M(\chi_b)$.
\begin{lmm}\label{lmb}
(1) \ Suppose $p \neq \pm{q},\pm{q^2},p^2 \neq -q^2$. Then $\Ind_{\widehat{\mathcal{H}}_1}^{\widehat{\mathcal{H}}}\rho_1^{b}$ and $\Ind_{\widehat{\mathcal{H}}_1}^{\widehat{\mathcal{H}}}\rho_2^{b}$ have 1- and 3-dimensional calibrated irreducible composition factors which are calibrated and presented by the following representation matrices; \\
(i) \ case $\Ind_{\widehat{\mathcal{H}}_1}^{\widehat{\mathcal{H}}}\rho_1^{b}$; 
\begin{itemize}
\item{$X_1=q^2p^{-1}, \ X_2=p^{-1}, \ T_1=q, \ T_2=-p^{-1}$. }
\item{$U_b^1$:
\begin{eqnarray*}
&{}& X_1=\left(
\begin{smallmatrix}
q^2p^{-1} & & \\
 & p & \\
 & & p
\end{smallmatrix}
\right),
 \ X_2=\left(
\begin{smallmatrix}
p & & \\
 & pq^{-2}& \\
 & & q^2p^{-1}
\end{smallmatrix}
\right), \\
&{}& T_1=\left(
\begin{smallmatrix}
\frac{q(q^2-1)}{q^2-p^2} & & -\frac{(p^2-1)(q^4-p^2)}{(q^2-p^2)^2}\\
 & q & \\
1 & & -\frac{p^2(q^2-1)}{(q^2-p^2)q} 
\end{smallmatrix}
\right), \ T_2=\left(
\begin{smallmatrix}
p & & \\
 & \frac{p(p^2-1)}{p^2-q^4} & 1 \\
 & -\frac{(p^2-q^2)(q^4-1)(p^2+q^2)}{(p^2-q^4)^2} & -\frac{(p^2-1)q^4}{p(p^2-q^4)} 
\end{smallmatrix}
\right). 
\end{eqnarray*}
}
\end{itemize}
(ii) \ case $\Ind_{\widehat{\mathcal{H}}_1}^{\widehat{\mathcal{H}}}\rho_2^{b}$; 
\begin{itemize}
\item{$X_1=pq^{-2}, \ X_2=p, \ T_1=-q^{-1}, \ T_2=-p$.} 
\item{$U_b^2$:
\begin{eqnarray*}
&{}&X_1=\left(
\begin{smallmatrix}
pq^{-2} & & \\
 & p^{-1} & \\
 & & p^{-1}
\end{smallmatrix}
\right),
 \ X_2=\left(
\begin{smallmatrix}
p^{-1} & & \\
 & pq^{-2}& \\
 & & q^2p^{-1}
\end{smallmatrix}
\right), \\
&{}& T_1=\left(
\begin{smallmatrix}
-\frac{p^2(q^2-1)}{q^2-p^2} & 1 & \\
-\frac{(p^2-1)(q^4-p^2)}{(q^2-p^2)^2} &  \frac{(q^2-1)q}{(q^2-p^2)} & \\
 & & -q^{-1}
\end{smallmatrix}
\right), \\
&{}& T_2=\left(
\begin{smallmatrix}
-p^{-1} & & \\
 & \frac{p(p^2-1)}{p^2-q^4} & 1 \\
 & -\frac{(p^2-q^2)(q^4-1)(p^2+q^2)}{(p^2-q^4)^2} & -\frac{(p^2-1)q^4}{p(p^2-q^4)} 
\end{smallmatrix}
\right). 
\end{eqnarray*}
}
\end{itemize}

\indent (2) \ Suppose $p=q$. Then they have 1-dimensional composition factor and 3-dimensional non-calibrated composition factor which are presented by the following representation matrices; \\
(i) \ case $\Ind_{\widehat{\mathcal{H}}_1}^{\widehat{\mathcal{H}}}\rho_1^{b}$;
\begin{itemize}
\item{$X_1=q, \ X_2=q^{-1}, \ T_1=q, \ T_2=-q^{-1}$.}
\item{$U_b^1$:
\begin{eqnarray*}
&{}& X_1=\left(
\begin{smallmatrix}
q & & \\
 & q & q^2  \\
 & & q
\end{smallmatrix}
\right),
X_2=\left(
\begin{smallmatrix}
q^{-1} & & \ \ \frac{1+2q^2}{q}\\
 & q & \ \ -q^2  \\
 & & \ \ q
\end{smallmatrix}
\right), \\
&{}& T_1=\left(
\begin{smallmatrix}
q & \ \  \frac{1+2q^2}{q^2} & \\
 & \ \ -q^{-1}  & \\
 & \ \ \frac{q^2-1}{q^2} & \ \ q
\end{smallmatrix}
\right),
T_2=\left(
\begin{smallmatrix}
-q^{-1} & & \ \  \frac{1+q^2}{q(q^2-1)}\\
1 & q & \ \  -\frac{1}{q^2-1} \\
 & & \ \ q
\end{smallmatrix}
\right). 
\end{eqnarray*}
}
\end{itemize}
(ii) \ case $\Ind_{\widehat{\mathcal{H}}_1}^{\widehat{\mathcal{H}}}\rho_2^{b}$; \\
\begin{itemize}
\item{$X_1=q^{-1}, \ X_2=q, \ T_1=-q^{-1}, \ T_2=q$.} 
\item{$U_b^2$:
\begin{eqnarray*}
&{}& X_1=\left(
\begin{smallmatrix}
q^{-1} & & -\frac{q^2-1}{q^3}\\
 & q^{-1} &   \\
 & & q^{-1}
\end{smallmatrix}
\right), \  
X_2=\left(
\begin{smallmatrix}
q^{-1} & & \ \ \frac{q^2-1}{q^3}\\
 & q & \ \ \frac{(q^2-1)(q^2+2)}{q}  \\
 & & \ \ q^{-1}
\end{smallmatrix}
\right), \\ 
&{}&T_1=\left(
\begin{smallmatrix}
q & & \\
q(2+q^2) & \ \ -q^{-1}  & \\
-q & & \ \ -q^{-1}
\end{smallmatrix}
\right), \  
T_2=\left(
\begin{smallmatrix}
-q^{-1} & \ \ q^{-1}& \ q\\
 & \ \ q & \ \ q(q^2+1) \\
 & & \ \ -q^{-1}
\end{smallmatrix}
\right).
\end{eqnarray*}
}
\end{itemize}
\indent (3) \ Suppose $p=q^2$. Then they have 1-dimensional composition factor and 3-dimensional non-calibrated composition factor which are presented by the following representation matrices; \\
(i) \ case $\Ind_{\widehat{\mathcal{H}}_1}^{\widehat{\mathcal{H}}}\rho_1^{b}$;
\begin{itemize}
\item{$X_1=1, \ X_2=q^{-2}, \ T_1=q, \ T_2=-q^{-2}$.}
\item{$U_b^1$:
\begin{eqnarray*}
&{}& X_1=\left(
\begin{smallmatrix}
1 & & \\
 & q^2 &   \\
 & & q^2
\end{smallmatrix}
\right), \  
X_2=\left(
\begin{smallmatrix}
q^2 & & \\
 & \ \ 1 & \ \ \frac{q^4-1}{q^2}  \\
 & & \ \ 1
\end{smallmatrix}
\right), \\ 
&{}&T_1=\left(
\begin{smallmatrix}
-q^{-1} & & \ \ -\frac{(q^2+1)^2}{q^2}\\
1 & \ \ q & \ \ \frac{q^2+1}{q}\\
 & & \ \ q
\end{smallmatrix}
\right), \  
T_2=\left(
\begin{smallmatrix}
q^2 & & \\
 & & 1 \\
 & 1 &  \ \ \frac{q^4-1}{q^2}
\end{smallmatrix}
\right).
\end{eqnarray*}
}
\end{itemize}
(ii) \ case  $\Ind_{\widehat{\mathcal{H}}_1}^{\widehat{\mathcal{H}}}\rho_2^{b}$: \\
\begin{itemize}
\item{$X_1=1, \ X_2=q^2, \ T_1=-q^{-1}, \ T_2=q^2$.}
\item{$U_b^2$:
\begin{eqnarray*}
&{}& X_1=\left(
\begin{smallmatrix}
q^{-2} & & \\
 & q^{-2} &   \\
 & & 1
\end{smallmatrix}
\right), \  
X_2=\left(
\begin{smallmatrix}
1 & \ \ \frac{q^4-1}{q^2}& \\
 &  \ \ 1 &   \\
 & & q^{-2}
\end{smallmatrix}
\right), \\ 
&{}&T_1=\left(
\begin{smallmatrix}
-q^{-1} &  \ \ \frac{q^2+1}{q} &  \ \ \frac{(q^2+1)^2}{q^2}\\
 &  \ \ -q^{-1} & \\
 &  \ \ 1&  \ \ q
\end{smallmatrix}
\right), \  
T_2=\left(
\begin{smallmatrix}
 &  \ \ 1 & \\
1 &  \ \ \frac{q^4-1}{q^2}&  \\
 & &  \ \ -q^{-2}
\end{smallmatrix}
\right).
\end{eqnarray*}
}
\end{itemize}
\end{lmm}
\end{exa}
\begin{exa}\label{exc}
We consider $M(\chi_c)$. Let $\rho_1^{c}$ and $\rho_2^{c}$ be the following 1-dimensional representations of $\widehat{\mathcal{H}}_2$;
\[
\begin{array}{|c||c|c|c|}
\hline
 & X_1 & X_2 & T_2 \\
\hline
\hline
\rho_1^{c} & -p^{-1} & p & p \\
\hline
\rho_2^{c} & -p^{-1} & -p^{-1} & -p^{-1}\\
\hline
\end{array}
\]
Since $\chi_c(\alpha_2^{\vee})=p^2$, we can apply the decompose lemma (Lemma \ref{dec}) to $M(\chi_c)$. 
\begin{lmm}\label{lmc}
(1) \ Suppose $p^2 \neq -q^{\pm{2}}$. $\Ind_{\widehat{\mathcal{H}}_2}^{\widehat{\mathcal{H}}}\rho_1^{c}$ and $\Ind_{\widehat{\mathcal{H}}_2}^{\widehat{\mathcal{H}}}\rho_2^{c}$ have two 2-dimensional irreducible calibrated composition factors which are presented by the following representation matricies;
\begin{eqnarray*}
&{}&composition \ factors \ of \ \Ind_{\widehat{\mathcal{H}}_2}^{\widehat{\mathcal{H}}}\rho_1^{c}; \\
&{}& \begin{array}{|c||c|c|c|c|}
\hline
 & X_1 & X_2 & T_1 & T_2 \\
\hline
\hline
U_c^1 & \left(
\begin{smallmatrix}
p& \\
&-p
\end{smallmatrix}
\right) & \left(
\begin{smallmatrix}
-p& \\
&p
\end{smallmatrix}
\right) & \left(
\begin{smallmatrix}
\frac{q^2-1}{2q}&\frac{(1+q^2)^2}{4q^2}\\
1&\frac{q^2-1}{2q}
\end{smallmatrix}
\right) & \left(
\begin{smallmatrix}
p& \\
&p
\end{smallmatrix}
\right) \\
\hline
U_c^3 & \left(
\begin{smallmatrix}
-p^{-1}& \\
&p
\end{smallmatrix}
\right) & \left(
\begin{smallmatrix}
p& \\
&-p^{-1}
\end{smallmatrix}
\right) & \left(
\begin{smallmatrix}
\frac{q^2-1}{(p^2+1)q}&\frac{(p^2+q^2)(1+p^2q^2)}{(p^2+1)^2q^2}\\
1&\frac{p^2(q^2-1)}{(p^2+1)q}
\end{smallmatrix}
\right) & \left(
\begin{smallmatrix}
p& \\
&-p^{-1}
\end{smallmatrix}
\right) \\
\hline
\end{array} \\
&{}&composition \ factors \ of \ \Ind_{\widehat{\mathcal{H}}_2}^{\widehat{\mathcal{H}}}\rho_2^{c}; \\
&{}& \begin{array}{|c||c|c|c|c|}
\hline
& X_1 & X_2 & T_1 & T_2 \\
\hline
\hline
U_c^2 & \left(
\begin{smallmatrix}
-p^{-1}& \\
&p^{-1}
\end{smallmatrix}
\right) & \left(
\begin{smallmatrix}
p^{-1}& \\
&-p^{-1}
\end{smallmatrix}
\right) & \left(
\begin{smallmatrix}
\frac{q^2-1}{2q}&\frac{(1+q^2)^2}{4q^2}\\
1&\frac{q^2-1}{2q}
\end{smallmatrix}
\right) & \left(
\begin{smallmatrix}
-p^{-1}& \\
&-p^{-1}
\end{smallmatrix}
\right) \\
\hline
U_c^4 & \left(
\begin{smallmatrix}
p^{-1}& \\
&-p
\end{smallmatrix}
\right) & \left(
\begin{smallmatrix}
-p& \\
&p^{-1}
\end{smallmatrix}
\right) & \left(
\begin{smallmatrix}
\frac{q^2-1}{(p^2+1)q}&\frac{(p^2+q^2)(1+p^2q^2)}{(p^2+1)^2q^2}\\
1&\frac{p^2(q^2-1)}{(p^2+1)q}
\end{smallmatrix}
\right) & \left(
\begin{smallmatrix}
p& \\
&-p^{-1}
\end{smallmatrix}
\right) \\

\hline
\end{array} \\
\end{eqnarray*}
\indent (2) \ Suppse $p^2=-q^2$. They have one 2-dimensional irreducible calibrated composition factor and two 1-dimensional composition factors. And their representation matrices are obtained by putting $p^2=-q^2$ in above matrices, since specialization lemma (Lemma \ref{spl}). More precisely, $U_c^1,U_c^2$ are irreducible, but $U_c^3,U_c^4$ have two 1-dimensional composition factors.
\end{lmm}
\end{exa}
\subsection{Classification Theorem}
By the preceding Examples and Lemmas, we obtain the following classification theorem. \\
\indent First, let us define the 1-dimensional representations of $\widehat{\mathcal{H}}_i$ in addition to the notation in the preceding Examples and Lemmas;
\begin{eqnarray*}
\begin{array}{|c||c|c|c|c|c|c|}
\hline
\widehat{\mathcal{H}}_1 & \rho_1^{d^{(1)}} & \rho_2^{d^{(1)}} & \rho_1^{d^{(2)}} & \rho_2^{d^{(2)}} & \rho_1^{g}(u) & \rho_2^{g}(u) \\
\hline
\hline
X_1 & q^2 & 1 &  q & q^{-1} & q^2u & u \\
\hline
X_2 & 1 & q^2 &  q^{-1} & q & u & q^2u \\
\hline
T_1 & q & -q^{-1} & q & -q^{-1} & q & -q^{-1} \\
\hline
\end{array}
\end{eqnarray*}
\begin{eqnarray*}
\begin{array}{|c||c|c|c|c|c|c|}
\hline
\widehat{\mathcal{H}}_2 & \rho_1^{d^{(3)}} & \rho_2^{d^{(3)}} & \rho_1^{d^{(4)}} & \rho_2^{d^{(4)}} & \rho_1^{f}(v) & \rho_2^{f}(v)\\
\hline
\hline
X_1 & p & p & 1 & 1 & pv & pv \\
\hline
X_2 & p & p^{-1} &  p & p^{-1} & p & p^{-1} \\
\hline
T_2 & p & -p^{-1} & p & -p^{-1} & p & -p^{-1} \\
\hline
\end{array}
\end{eqnarray*}
\begin{thm}\label{MT}
Suppose that $p$ and $q$ are not a root of unity. The finite-dimensional irreducible representations of type $B_2$ with unequal parameters are given by the following lists depending on the relation of parameters. \\
(0) \ The principal series representations $M(\chi)$, where $\chi \neq \pm{\chi_a},\pm{\chi_b},\chi_c,\pm{\chi_d^{(j)}} \ (1 \le j \le 5),\pm{\chi_f(v)},\chi_g(u)$ and their W-orbits, are irreducible. \\
(1) \ For any $p,q$, there are eight 1-dimensional (irreducible) representations defined by 
\[
\begin{array}{|c||c|c|c|c|c|c|c|c|}
\hline
X_1 & q^2p & q^{-2}p^{-1} & q^2p^{-1} & q^{-2}p & -q^2p & -q^{-2}p^{-1} & -q^2p^{-1} & -q^{-2}p \\
\hline
X_2 & p & p^{-1} & p^{-1} & p & -p & -p^{-1} & -p^{-1} & -p \\
\hline
T_1 & q & -q^{-1} & q & -q^{-1} & q & -q^{-1} & q & -q^{-1} \\
\hline
T_2 & p & -p^{-1} & -p^{-1} & p & p & -p^{-1} & -p^{-1} & p \\
\hline
\end{array}
\]
(2) \ For any $p,q$,
\begin{eqnarray*}
&{}&  \ \Ind_{\widehat{\mathcal{H}}_2}^{\widehat{\mathcal{H}}}\rho_1^{f}(v), \ \Ind_{\widehat{\mathcal{H}}_2}^{\widehat{\mathcal{H}}}\rho_2^{f}(v),\ \Ind_{\widehat{\mathcal{H}}_2}^{\widehat{\mathcal{H}}}(-\rho_1^{f}(v)), \ \Ind_{\widehat{\mathcal{H}}_2}^{\widehat{\mathcal{H}}}(-\rho_2^{f}(v))\\ 
&{}& \quad \quad \text{with} \ v \neq \pm{p^{-2}},\pm{p^{-1}},\pm{1},q^{\pm{2}},q^{\pm{2}}p^{-2}\\
&{}& \ \Ind_{\widehat{\mathcal{H}}_1}^{\widehat{\mathcal{H}}}\rho_1^{g}(u), \ \ \Ind_{\widehat{\mathcal{H}}_1}^{\widehat{\mathcal{H}}}\rho_2^{g}(u) \ \text{with} \ u \neq \pm{p^{\pm{1}}},\pm{1},\pm{q^{-2}},\pm{q^{-1}},\pm{q^{-2}p^{\pm{1}}}
\end{eqnarray*}
are 4-dimensional one parameter families of irreducible representations and calibrated. They are not isomorphic to each other. \\
(3) \ When $p,q$ are generic i.e. $p \neq \pm{q^{\pm{2}}},\pm{q^{\pm{1}}}$ and $p^2 \neq -q^{\pm{2}}$, the remaining finite-dimensional irreducible representations are the following; \\
\indent (I) \ $U_c^i \ (1 \le i \le 4)$ which are 2-dimensional and calibrated. \\
\indent (II) \ $U_a^i,U_b^i,U_{-a}^i,U_{-b}^i \ (i=1,2)$ which are 3-dimensional and calibrated. \\
\indent (III) \ 
\begin{eqnarray*}
&{}&\Ind_{\widehat{\mathcal{H}}_1}^{\widehat{\mathcal{H}}}\rho_j^{d^{(i)}},\Ind_{\widehat{\mathcal{H}}_1}^{\widehat{\mathcal{H}}}(-\rho_j^{d^{(i)}}) \ (j=1,2,i=1,2), \\
&{}&\Ind_{\widehat{\mathcal{H}}_2}^{\widehat{\mathcal{H}}}\rho_j^{d^{(i)}},\Ind_{\widehat{\mathcal{H}}_2}^{\widehat{\mathcal{H}}}(-\rho_j^{d^{(i)}}) \ (j=1,2,i=3,4,5) 
\end{eqnarray*}
which are 4-dimensional and non-calibrated. \\
(4) \ When $p=q^2$, the remaining finite-dimensional irreducible representations are the following; \\
\indent (I) \ $U_c^i \ (1 \le i \le 4)$ which are 2-dimensional and calibrated. \\
\indent (II) \ $U_a^i,U_{-a}^i, \ (i=1,2)$ which are 3-dimensional and calibrated. \\
\indent (III) \ $U_b^i,U_{-b}^i, \ (i=1,2)$ which are 3-dimensional and non-calibrated. \\
\indent (IV) \ 
\begin{eqnarray*}
&{}&\Ind_{\widehat{\mathcal{H}}_1}^{\widehat{\mathcal{H}}}\rho_j^{d^{(i)}},\Ind_{\widehat{\mathcal{H}}_1}^{\widehat{\mathcal{H}}}(-\rho_j^{d^{(i)}}) \ (j=1,2,i=2), \\
&{}&\Ind_{\widehat{\mathcal{H}}_2}^{\widehat{\mathcal{H}}}\rho_j^{d^{(i)}},\Ind_{\widehat{\mathcal{H}}_2}^{\widehat{\mathcal{H}}}(-\rho_j^{d^{(i)}}) \ (j=1,2,i=3,5) 
\end{eqnarray*}
which are 4-dimensional and non-calibrated. \\
(5) \ When $p=q$, the remaining finite-dimensional irreducible representations are the following; \\
\indent (I) \ $U_c^i \ (1 \le i \le 4)$ which are 2-dimensional and calibrated. \\
\indent (II) \ $U_a^i,U_{-a}^i, \ (i=1,2)$ which are 3-dimensional and calibrated. \\
\indent (III) \ $U_b^i,U_{-b}^i, \ (i=1,2)$ which are 3-dimensional and non-calibrated. \\
\indent (IV) \ 
\begin{eqnarray*}
&{}&\Ind_{\widehat{\mathcal{H}}_1}^{\widehat{\mathcal{H}}}\rho_j^{d^{(i)}},\Ind_{\widehat{\mathcal{H}}_1}^{\widehat{\mathcal{H}}}(-\rho_j^{d^{(i)}}) \ (j=1,2,i=1), \\
&{}&\Ind_{\widehat{\mathcal{H}}_2}^{\widehat{\mathcal{H}}}\rho_j^{d^{(i)}},\Ind_{\widehat{\mathcal{H}}_2}^{\widehat{\mathcal{H}}}(-\rho_j^{d^{(i)}}) \ (j=1,2,i=4,5) 
\end{eqnarray*}
which are 4-dimensional and non-calibrated. \\

(6) \ When $p^2=-q^2$, the remaining finite-dimensional irreducible representations are the following; \\
\indent (I) \ $U_c^i \ (i=1,2)$ which are 2-dimensional and calibrated. \\
\indent (II) \ $U_a^i,U_{-a}^i, \ (1 \le i \le 2)$ which are 3-dimensional and calibrated. \\
\indent (III) \ 
\begin{eqnarray*}
&{}&\Ind_{\widehat{\mathcal{H}}_1}^{\widehat{\mathcal{H}}}\rho_j^{d^{(i)}},\Ind_{\widehat{\mathcal{H}}_1}^{\widehat{\mathcal{H}}}(-\rho_j^{d^{(i)}}) \ (j=1,2,i=1,2), \\
&{}&\Ind_{\widehat{\mathcal{H}}_2}^{\widehat{\mathcal{H}}}\rho_j^{d^{(i)}},\Ind_{\widehat{\mathcal{H}}_2}^{\widehat{\mathcal{H}}}(-\rho_j^{d^{(i)}}) \ (j=1,2,i=3,4,5) 
\end{eqnarray*}
which are 4-dimensional and non-calibrated. \\
(7) \ Using the following automorphisms of $\mathcal{\widehat{H}}$
\[
X_1 \mapsto X_1,X_2 \mapsto X_2,T_1 \mapsto T_1,T_2 \mapsto -T_2, q \mapsto q, p \mapsto \mp{p}^{\pm{1}}
\]
the cases of $p=\pm{q^{-2}},-q^2$ reduces the case (4). Similarly, the cases of $p=\pm{q^{-1}},-q$ reduces the case (5). The case of $p^2=-q^{-2}$ also reduces the case (6).\\
\end{thm}
\begin{note}
In \cite{R2}, Ram dealt equal parameter case i.e. $p=q$ case. However the case $\chi_d^{(5)}$ does not appear in his list. Also he did not explicitly list $-\chi_a,-\chi_b,-\chi_d^{(1)},-\chi_d^{(4)},-\chi_d^{(5)}$ and $-\chi_f$. 
\end{note}
\section{Tables of irreducible representations}
We will summarize about the dimension of composition factors and their calibratability. Note that we will omit the principal series representation $M(-\chi)$ and their composition factors in the following tables. 
\newpage
{\footnotesize
\subsection{$p,q$ generic case (i.e. $p \neq \pm{q}^{\pm{1}},\pm{q}^{\pm{2}}$ and $p^2 \neq -q^{\pm{2}}$)}
\begin{eqnarray*}
\begin{array}{|c||c|c|c|c|c|}
\hline
 & \chi(X_1) & \chi(X_2) & P(\chi) & \text{dim} & \text{calibrated?} \\
\hline
\hline
\chi_a & q^2p & p & \{\alpha_1,\alpha_2\} & 1 & \bigcirc \\
    &      &   &                       & 3 & \bigcirc \\
    &      &   &                       & 3 & \bigcirc \\
    &      &   &                       & 1 & \bigcirc \\
\hline
\chi_b & q^2p^{-1} & p^{-1} & \{\alpha_1,\alpha_2\} & 1 & \bigcirc  \\
    &      &   &                       & 3 & \bigcirc \\
    &      &   &                       & 3 & \bigcirc \\
    &      &   &                       & 1 & \bigcirc \\
\hline
\chi_c & -p^{-1} & p & \{\alpha_2,2\alpha_1+\alpha_2\} & 2 & \bigcirc \\
    &      &   &                       & 2 & \bigcirc \\
    &      &   &                       & 2 & \bigcirc \\
    &      &   &                       & 2 & \bigcirc \\
\hline
\chi_d^{(1)} & q^2 & 1 & \{\alpha_1,\alpha_1+\alpha_2\} & 4 & \times \\
    &      &   &                       & 4 & \times \\
\hline
\chi_d^{(2)} & q & q^{-1}& \{\alpha_1\}   & 4 & \times \\
    &      &   &                       & 4 & \times \\
\hline
\chi_d^{(3)} & p & p & \{\alpha_2,2\alpha_1+\alpha_2\} & 4 & \times \\
    &      &   &                       & 4 & \times \\
\hline
\chi_d^{(4)} & 1 & p& \{\alpha_2\}               & 4 & \times \\
    &      &   &                       & 4 & \times \\
\hline
\chi_d^{(5)} & -1 & p & \{\alpha_2\} & 4 & \times \\
    &      &   &                       & 4 & \times \\
\hline
\chi_f(v) & pv & p & \{\alpha_2\}               & 4 & \bigcirc \\
    &      &   &                       & 4 & \bigcirc \\
\hline
\chi_g(u) & q^2u & u & \{\alpha_1\}               & 4 & \bigcirc \\
    &      &   &                & 4 & \bigcirc \\
\hline
\end{array}
\end{eqnarray*}
\subsection{$p=q$ case; equal parameter case}
\begin{eqnarray*}
\begin{array}{|c||c|c|c|c|c|}
\hline
 & \chi(X_1) & \chi(X_2) & P(\chi) & \text{dim} & \text{calibrated?} \\
\hline
\hline
\chi_a & q^3 & q & \{\alpha_1,\alpha_2\} & 1 & \bigcirc \\
    &      &   &                       & 3 & \bigcirc \\
    &      &   &                       & 3 & \bigcirc \\
    &      &   &                       & 1 & \bigcirc \\
\hline
\chi_b & q & q^{-1} & \{\alpha_1,\alpha_2,2\alpha_1+\alpha_2\} & 1 & \bigcirc  \\
    &      &   &                       & 3 & \times \\
    &      &   &                       & 3 & \times \\
    &      &   &                       & 1 & \bigcirc \\
\hline
\chi_c & -q^{-1} & q & \{\alpha_2,2\alpha_1+\alpha_2\} & 2 & \bigcirc \\
    &      &   &                       & 2 & \bigcirc \\
    &      &   &                       & 2 & \bigcirc \\
    &      &   &                       & 2 & \bigcirc \\
\hline
\chi_d^{(1)} & q^2 & 1 & \{\alpha_1,\alpha_1+\alpha_2\} & 4 & \times \\
    &      &   &                       & 4 & \times \\
\hline
\chi_d^{(4)} & 1 & q& \{\alpha_2\}               & 4 & \times \\
    &      &   &                       & 4 & \times \\
\hline
\chi_d^{(5)} & -1 & p & \{\alpha_2\} & 4 & \times \\
    &      &   &                       & 4 & \times \\
\hline
\chi_f(v) & qv & q & \{\alpha_2\}               & 4 & \bigcirc \\
    &      &   &                       & 4 & \bigcirc \\
\hline
\chi_g(u) & q^2u & u & \{\alpha_1\}               & 4 & \bigcirc \\
    &      &   &                & 4 & \bigcirc \\
\hline
\end{array}
\end{eqnarray*}
\newpage
\subsection{$p=q^2$ case}
\begin{eqnarray*}
\begin{array}{|c||c|c|c|c|c|}
\hline
 & \chi(X_1) & \chi(X_2) & P(\chi) & \text{dim} & \text{calibrated?} \\
\hline
\hline
\chi_a & q^4 & q^2 & \{\alpha_1,\alpha_2\} & 1 & \bigcirc \\
    &      &   &                       & 3 & \bigcirc \\
    &      &   &                       & 3 & \bigcirc \\
    &      &   &                       & 1 & \bigcirc \\
\hline
\chi_b & 1 & q^{-2} & \{\alpha_1,\alpha_2,\alpha_1+\alpha_2\} & 1 & \bigcirc  \\
    &      &   &                       & 3 & \times \\
    &      &   &                       & 3 & \times \\
    &      &   &                       & 1 & \bigcirc \\
\hline
\chi_c & -q^{-2} & q^2 & \{\alpha_2,2\alpha_1+\alpha_2\} & 2 & \bigcirc \\
    &      &   &                       & 2 & \bigcirc \\
    &      &   &                       & 2 & \bigcirc \\
    &      &   &                       & 2 & \bigcirc \\
\hline
\chi_d^{(2)} & q & q^{-1}& \{\alpha_1\}   & 4 & \times \\
    &      &   &                       & 4 & \times \\
\hline
\chi_d^{(3)} & q^2 & q^2 & \{\alpha_2,2\alpha_1+\alpha_2\} & 4 & \times \\
    &      &   &                       & 4 & \times \\
\hline
\chi_d^{(5)} & -1 & q^2 & \{\alpha_2\} & 4 & \times \\
    &      &   &                       & 4 & \times \\
\hline
\chi_f(v) & q^2v & q^2 & \{\alpha_2\}               & 4 & \bigcirc \\
    &      &   &                       & 4 & \bigcirc \\
\hline
\chi_g(u) & q^2u & u & \{\alpha_1\}               & 4 & \bigcirc \\
    &      &   &                & 4 & \bigcirc \\
\hline
\end{array}
\end{eqnarray*}
\subsection{$p^2=-q^2$ case}
\begin{eqnarray*}
\begin{array}{|c||c|c|c|c|c|}
\hline
 & \chi(X_1) & \chi(X_2) & P(\chi) & \text{dim} & \text{calibrated?} \\
\hline
\hline
\chi_a & -p^3 & p & \{\alpha_1,\alpha_2\} & 1 & \bigcirc \\
    &      &   &                       & 3 & \bigcirc \\
    &      &   &                       & 3 & \bigcirc \\
    &      &   &                       & 1 & \bigcirc \\
\hline
\chi_c & -p^{-1} & p & \{\alpha_1,\alpha_2,2\alpha_1+\alpha_2\} & 1 & \bigcirc \\
    &      &   &                       & 1 & \bigcirc \\
    &      &   &                       & 1 & \bigcirc \\
    &      &   &                       & 1 & \bigcirc \\
    &      &   &                       & 2 & \bigcirc \\
    &      &   &                       & 2 & \bigcirc \\
\hline
\chi_d^{(1)} & -p^2 & 1 & \{\alpha_1,\alpha_1+\alpha_2\} & 4 & \times \\
    &      &   &                       & 4 & \times \\
\hline
\chi_d^{(2)} & \pm{p\sqrt{-1}} & \pm{p}\sqrt{-1}& \{\alpha_1\}   & 4 & \times \\
    &      &   &                       & 4 & \times \\
\hline
\chi_d^{(3)} & p & p & \{\alpha_2,2\alpha_1+\alpha_2\} & 4 & \times \\
    &      &   &                       & 4 & \times \\
\hline
\chi_d^{(4)} & 1 & p& \{\alpha_2\}               & 4 & \times \\
    &      &   &                       & 4 & \times \\
\hline
\chi_d^{(5)} & -1 & p & \{\alpha_2\} & 4 & \times \\
    &      &   &                       & 4 & \times \\
\hline
\chi_f(v) & pv & p & \{\alpha_2\}               & 4 & \bigcirc \\
    &      &   &                       & 4 & \bigcirc \\
\hline
\chi_g(u) & -p^2u & u & \{\alpha_1\}               & 4 & \bigcirc \\
    &      &   &                & 4 & \bigcirc \\
\hline
\end{array}
\end{eqnarray*}
}

\end{document}